\documentclass [11pt]{amsart}
\usepackage {amsmath, amssymb, amscd, graphicx, color, url,epsf}
\usepackage[all,knot,arc]{xy}

\newtheorem {theorem}{Theorem}[section]
\newtheorem {definition}[theorem]{Definition}
\newtheorem {lemma}[theorem]{Lemma}
\newtheorem {proposition}[theorem]{Proposition}
\newtheorem {corollary}[theorem]{Corollary}

\newcommand\Z{\mathbb{Z}}
\newcommand\R{\mathbb{R}}
\newcommand\Ta{\mathbb{T}_\alpha}
\newcommand\Tb{\mathbb{T}_\beta}

\def\Sym{\mathrm{Sym}}
\def\Filt{\mathcal F}
\def\HFa{\widehat {HF}}
\def\HFm{{HF}^-}
\def\CFm{{CF}^-}
\def\CFKm{{CFK}^-}
\def\CFKa{\widehat{CFK}}
\def\CFKinf{{CFK}^\infty}
\def\CFKp{{CFK}^+}
\def\CFa{\widehat{CF}}

\def\HFLa{{\widehat{\mathrm{HFL}}}}
\def\HFKa{{\widehat{\mathrm{HFK}}}}
\def\HFLa{{\widehat{\mathrm{HFL}}}}
\def\H{\mathbb H}
\def\OneHalf{\frac{1}{2}}

\def\rk {{\operatorname{rank}}}

\def\cald{\mathcal D}
\def\CFKm{{CFK}^-}
\def\CFLm{{CFL}^-}

\def\cm{\cdot}
\def\x{\mathbf x}
\def\y{\mathbf y}
\def\z{\mathbf z}
\def\Torus{\mathcal T}
\def\ModFlow{\mathcal M}
\newcommand\orL{\vec{L}}
\newcommand\Mas{\mu}
\def\Field{\mathbb F}

\begin{document}

\title
[A combinatorial description of knot Floer homology]{A combinatorial
description of knot Floer homology}

\author [Ciprian Manolescu]{Ciprian Manolescu}
\thanks {CM was supported by a Clay Research Fellowship.}
\address {Department of Mathematics, Columbia University\\ New York, NY 10027}
\email {cm@math.columbia.edu}

\author [Peter Ozsv\'ath]{Peter Ozsv\'ath}
\thanks {PSO was supported by NSF grant number DMS-0505811 and FRG-0244663}
\address {Department of Mathematics, Columbia University\\ New York, NY 10027}
\email {petero@math.columbia.edu}

\author [Sucharit Sarkar]{Sucharit Sarkar}
\address {Department of Mathematics, Princeton Univeristy\\ Princeton, NJ 08544}
\email {sucharit@math.princeton.edu}

\begin {abstract} 
Given a grid presentation of a knot (or link) $K$ in the three-sphere,
we describe a Heegaard diagram for the knot complement in which the
Heegaard surface is a torus and all elementary domains are
squares. Using this diagram, we obtain a purely combinatorial
description of the knot Floer homology of $K$.
\end {abstract}

\maketitle
\section{Introduction}
\label{sec:Introduction}

Heegaard Floer homology~\cite{HolDisk} is an invariant for
three-manifolds, defined using holomorphic disks and Heegaard
diagrams.  In~\cite{Knots} and~\cite{RasmussenThesis}, this
construction is extended to give an invariant, {\em knot Floer
  homology} $\HFKa$, for null-homologous knots in a closed, oriented
three-manifold. This construction is further generalized
in~\cite{Links} to the case of oriented links.  The
definition of all these invariants involves counts of holomorphic disks
in the symmetric product of a Riemann surface, which makes them rather
challenging to calculate.

In its most basic form, knot Floer homology is an invariant for knots
$K\subset S^3$, $\HFKa(K)$, which is a finite-dimensional bi-graded
vector space over $\Field={\mathbb Z}/2{\mathbb Z}$, i.e.
$$\HFKa(K)=\bigoplus_{m,s} \HFKa_m(K,s).$$ This invariant is related to the
symmetrized Alexander polynomial $\Delta_K(T)$ by the formula
\begin{equation}
  \label{eq:RelateWithAlexander}
\Delta_K(T)=\sum_{m,s} (-1)^m \rk~\HFKa_m(K,s) \cm T^s
\end{equation}
(cf.~\cite{Knots},
\cite{RasmussenThesis}). The topological significance of this invariant
is illustrated by 
the result that
$$g(K)=\max\{s\in\Z\big| \HFKa_*(K,s)\neq 0\},$$
where here $g(K)$
denotes the Seifert genus of $K$
(cf.~\cite{GenusBounds}), and also the fact that
$\HFKa_*(K,g(K))$ has rank one if and only if $K$ is fibered
(\cite{Ghiggini}
in the case where $g(K)=1$ and \cite{YiNiFibered} in general).
The invariant is defined as a version of Lagrangian Floer
homology~\cite{FloerLag} in a suitable symmetric product of a Heegaard surface.

Our aim here is to give a purely combinatorial presentation of knot
Floer homology with coefficients in $\Field$ for knots in the
three-sphere. Our description can be extended to describe link Floer
homology, and also it can be extended to describe the ``full knot
filtration'' (and in particular the concordance invariant
$\tau$~\cite{4BallGenus}).  However, in the interest of exposition, we
limit ourselves in the introduction to the case of knot Floer
homology, referring the interested reader to
Section~\ref{sec:Conclusion} for more general cases.

To explain our combinatorial description, it will be useful to have
the following notions.

A {\em planar grid diagram} ${\widetilde\Gamma}$ consists of a square
grid on the the plane with $n \times n$ cells, together with a
collection of black and white dots on it, arranged so that:
\begin{itemize}
\item {every row contains exactly one black dot and one white dot;}
\item {every column contains exactly one black dot and one white dot;}
\item {no cell contains more than one dot.}
\end {itemize}
The number $n$ is called the {\em grid number} of ${\widetilde\Gamma}.$

\begin{figure}
\begin{center}
\mbox{\vbox{\epsfbox{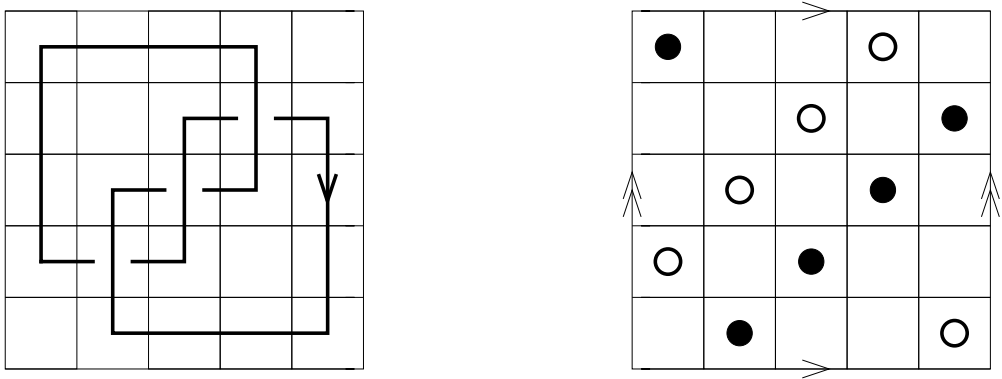}}}
\end{center}
\caption {{\bf Grid diagram for the trefoil.} We have pictured here a 
grid diagram for the trefoil, with projection indicated on the left.
To pass from a planar to a toroidal grid diagram, we make the
identifications suggested by the arrows.}
\label{Figure:tref}
\end{figure}

Given a planar grid diagram ${\widetilde\Gamma}$, we can place it in a
standard position on the plane as follows: the bottom left corner is
at the origin, each cell is a square of edge length one, and every dot
is in the middle of the respective cell. We then construct a planar
knot projection by drawing horizontal segments from the white to the
black dot in each row, and vertical segments from the black to the
white dot in each column. At every intersection point, we let the
horizontal segment be the underpass and the vertical one the
overpass. This produces a planar diagram for an oriented link $\orL$
in $S^3$. We say that $\orL$ has a grid presentation given by
${\widetilde\Gamma}.$ Figure~\ref{Figure:tref} shows a grid
presentation of the trefoil, with $n=5.$

It is easy to see that every knot (or link) in the three-sphere can be
presented by a planar grid diagram. In fact, grid presentations are
equivalent to the arc presentations of knots, which first appeared in
\cite{Brunn}, the square bridge positions of knots of~\cite{Lyon}, and
also to Legendrian realizations of knots, cf.~\cite{Matsuda}; they
have enjoyed a considerable amount of attention over the years, see
also~\cite{Cromwell}, \cite{Dynnikov}.  The
minimum number $n$ for which a knot $K\subset S^3$ admits a grid
presentation of grid number $n$ is called the {\em arc index} of $K$.

We find it convenient to transfer our planar grid diagrams to the
torus $\Torus$ obtained by gluing the topmost segment to the bottom-most one,
and the leftmost segment to the rightmost one. In the torus, our
horizontal and vertical arcs become horizontal and vertical circles.
The torus inherits its orientation from the plane.
We call the resulting object ${\Gamma}$ a {\em toroidal grid diagram},
or simply a grid diagram, for $K$.  

Given a toroidal grid diagram, we associate to it a chain complex
$\bigl({ C}(\Gamma), {\partial} \bigr)$ as follows. The generators $X$ of
${ C}(\Gamma)$ are indexed by one-to-one correspondences between the
horizontal and vertical circles. More geometrically, we can think of
these as $n$-tuples of intersection points $\x$ between the horizontal
and vertical circles, with the property that no intersection point
appears on more than one horizontal (or vertical) circle.

We now define functions $A\colon X \longrightarrow \Z$ and $M\colon X
\longrightarrow \Z$ (the Alexander and Maslov gradings) as follows.

\begin{figure}
\begin{center}
\mbox{\vbox{\epsfbox{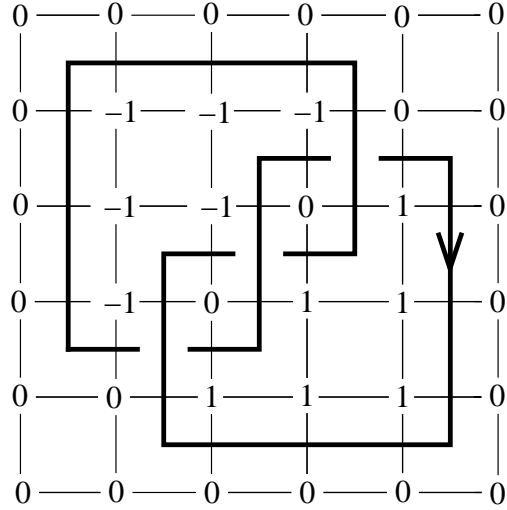}}}
\end{center}
\caption {{\bf The function $a$.} Over every lattice point $p$ from 
Figure~\ref{Figure:tref}, we marked minus the winding number of the 
knot projection around $p.$}
\label{Figure:Alexander}
\end{figure}

Let us define a function $a$ on lattice points $p$ to be minus one times 
the winding number of the knot projection around $p$. (This is shown in 
Figure~\ref{Figure:Alexander} for our trefoil example.) Each black or 
white dot in the diagram lies in a square. We thus obtain $2n$ 
distinguished squares, and each of them has four corners. We denote the 
resulting collection of corners $\{c_{i,j}\}, \ i\in\{1,...,2n\},\ 
j\in\{1,...,4\}.$ We set 
\begin {equation}
\label{eq:Alexanders}
A(\x)=\sum_{p\in \x} a(p)- \frac{1}{8} \Bigl(\sum_{i,j} 
a(c_{i,j})\Bigr) - \frac{n-1}{2}.
\end {equation}

Next, given a pair of generators $\x$ and $\y$, and an embedded rectangle 
$r$ in $\Torus$ whose edges are arcs in the horizontal and vertical 
circles, we say that $r$ connects $\x$ to $\y$ if $\x$ and $\y$ agree 
along all but two horizontal circles, if all four corners of $r$ are 
intersection points in $\x\cup\y$, and indeed, if we traverse each 
horizontal boundary components of $r$ in the direction dictated by the 
orientation that $r$ inherits from $\Torus$, then the arc is oriented so 
as to go from a point in $\x$ to the point in $\y$. Let $R_{\x,\y}$ denote 
the collection of rectangles connecting $\x$ to $\y$.

\begin{figure}
\begin{center}
\mbox{\vbox{\epsfbox{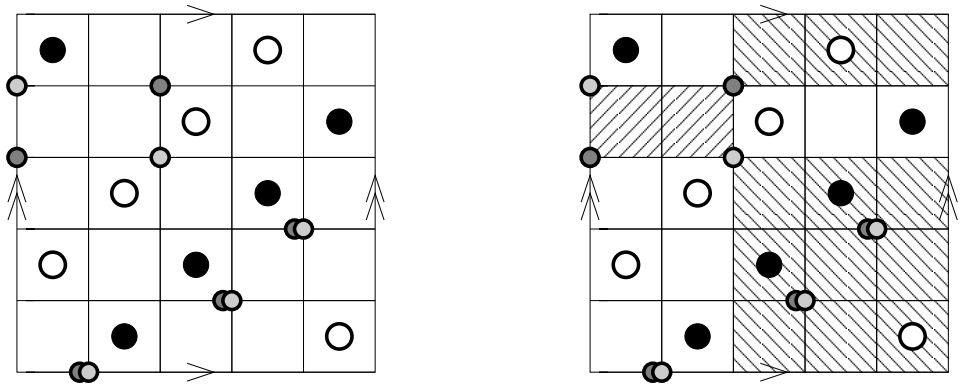}}}
\end{center}
\caption {{\bf Rectangles.} At the left, we have indicated two generators
$\x$ and $\y$ in $X$ for the grid diagram of the trefoil considered
earlier. The generators $\x$ and $\y$ are represented by the
collections of (smaller) shaded dots centered on the intersection
points of the grid, with $\x$ represented by the more darkly shaded
circles and $\y$ represented by the more lightly shaded ones. Note
that three dots in $\x$ occupy the same locations on the grid as
$\y$-dots, while two do not. At the right, we have have indicated the
two rectangles in $R_{\x,\y}$, which are shaded by (the two types of)
diagonal hatchings. One of these rectangles $r$ has
$P_{\x}(r)+P_{\y}(r)=1$ and $W(r)=B(r)=0$ (and hence it represents a
non-trivial differential from $\x$ to $\y$), while the other rectangle
$r'$ has $P_{\x}(r')+P_{\y}(r')=5$ and $W(r')=B(r')=2$.}
\label{Figure:rect}
\end{figure}

It is easy to see that if $\x,\y\in X$, and if $\x$ and $\y$ differ
along exactly two horizontal circles, then there are exactly two
rectangles in $R_{\x,\y}$; otherwise $R_{\x,\y}=\emptyset$, cf.
Figure~\ref{Figure:rect}.

\begin{figure}
\begin{center}
\mbox{\vbox{\epsfbox{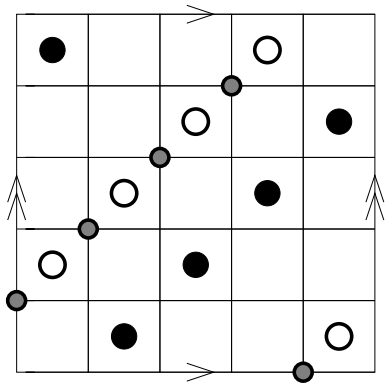}}}
\end{center}
\caption {{\bf The generator $\x_0$.} We have illustrated here
the generator with Maslov grading equal to $1-n$.}
\label{fig:MaslovZero}
\end{figure}

Given $\x,\y\in X$, it is easy to find an oriented, null-homologous
curve $\gamma_{\x,\y}$ composed of horizontal and vertical arcs, where
each horizontal arc goes from a point in $\x$ to a point in $\y$ (and
hence each vertical arc goes from a point in $\y$ to a point in $\x$).
Now, suppose that $D$ is a two-chain whose boundary is a collection of
horizontal and vertical arcs, and $\x\in X$. We let $W(D)$ and $B(D)$
denote the number of white and black dots in $D$ respectively.
Moreover, near each intersection point $x$ of the horizontal and
vertical circles, $D$ has four local multiplicities. We define the
local multiplicity of $D$ at $x$, $p_{x}(D)$, to be the average of
these four local multiplicities. Moreover, given $\x\in X$, let
$$P_{\x}(D)= \sum_{x\in \x} p_{x}(D).$$ Now, $M$ is uniquely
characterized up to an additive constant by the property that for each
$\x,\y\in X$, 
\begin{equation}
\label{eq:RelativeMaslovGrading}
M(\x)-M(\y)=P_{\x}(D)+P_{\y}(D)-2\cm W(D),
\end{equation}
where here $D$ is some two-chain whose boundary is $\gamma_{\x,\y}$.
(Observe that we have displayed here a simple special case of
Lipshitz's formula for the Maslov index of a holomorphic disk in the
symmetric product, cf.~\cite{LipshitzCyl}.) Note that the
right-hand-side is independent of the choice of $D$, as follows. Let
$\{A_i\}_{i=1}^n$ and $\{B_i\}_{i=1}^n$ be the annuli given by
$\Sigma-\alpha_1-...-\alpha_n$ and $\beta_1-...-\beta_n$ respectively.
Note that any two choices of $D$ and $D'$ connecting $\x$ to $\y$
differ by adding or subtracting a finite number of annuli $A_i$ and
$B_j$.  But for each such annulus $A$, $P_\x(A)=1$, $P_\y(A)=1$, and
$W(A)=1$, and hence they do not change the right-hand-side.  Moreover,
the additive indeterminacy in $M$ is removed by the following
convention.  Consider the generator $\x_0$ which occupies the lower
left-hand corner of each square which contains a white dot,
cf. Figure~\ref{fig:MaslovZero}. We declare that $M(\x_0)=1-n$.

Consider ${ C}(\Gamma)$, the $\mathbb{F}$-vector space generated by elements of
$X$. We define a differential
$$\partial \colon { C}(\Gamma) \longrightarrow { C}(\Gamma)$$
by the formula
$${\partial} \x = \sum_{\y\in X} \sum_{r\in R_{\x,\y}} \left\{\begin{array}{ll}
1 & {\text{if $P_{\x}(r)+P_{\y}(r)=1$ and $W(r)=B(r)=0$}} \\
0 & {\text{otherwise}}
\end{array}
\right\} \cm \y.$$
The condition that $P_{\x}(r)+P_{\y}(r)=1$ and $W(r)=B(r)=0$
is, of course, equivalent to the condition that the interior of the rectangle $r$
contains no black points, white points, or points amongst the $\x$ and $\y$.

It is easy to see that ${\partial}$ drops Maslov grading by one and
preserves Alexander grading. It is also elementary to verify that
$\partial^2=0$.  Thus, we can take the homology of this
complex to obtain a bigraded vector space over $\Field$.

Let $V$ be the two-dimensional bigraded vector space spanned
by one generator in bigrading $(-1,-1)$ and another in bigrading $(0,0)$.

We can now state the following:

\begin{theorem}
  \label{thm:KnotFloerHomology} Fix a grid presentation $\Gamma$ of a
  knot $K$, with grid number $n$. Then, the homology of the above
  chain complex $H_*({ C}(\Gamma),{\partial})$ is isomorphic to
  the bigraded group $\HFKa(K)\otimes V^{\otimes (n-1)}$.
\end{theorem}

The key point of the above theorem is to find a suitable Heegaard
diagram for $S^3$ compatible with the knot $K$. Indeed, the diagram we
use has genus one, with Heegaard torus $\Torus$, and $K$ is
represented as a collection of horizontal and vertical arcs. This
Heegaard diagram has the property that the knot pierces $\Torus$ in 
several pairs of points, and the very interesting property that the
complement in $\Torus$ of the attaching circles is a collection of
squares. In this case, properties of the Maslov index ensure that the
only holomorphic disks are rectangles. The chain complex ${ C}(\Gamma)$
we have described above, then, agrees with the Heegaard Floer complex
for this diagram.

There are several other variants of
Theorem~\ref{thm:KnotFloerHomology}. There is, for example, a version
which calculates $\HFKa(K)$ directly (though, of course, it is
uniquely determined by the above result), except there one needs to
consider a variant of the above the chain complex 
defined over a suitable polynomial
algebra. 

We also consider in this paper several other versions of
Theorem~\ref{thm:KnotFloerHomology}. We discuss how to calculate the
other variants of knot Floer homology, and also a variant for
links.

This paper is organized as follows. In Section~\ref{sec:Basepoints},
we describe the construction of link Floer homology using Heegaard
diagrams with the property that the link crosses the Heegaard surface
in many points. This construction is then identified with the usual
construction using methods from~\cite{Links}. In
Section~\ref{sec:Conclusion}, we identify the chain complex ${C}(\Gamma)$
with the link Floer homology complex using the toroidal grid diagram of $L$,
interpreted as a Heegaard diagram for $L$, and state some more general
consequences. Finally, in Section~\ref{sec:Examples}, we describe some
simple examples to illustrate our results.

\vskip.3cm
\noindent{\textbf{Further remarks.}}
Whereas the constructions in this paper give a purely combinatorial
chain complex for knot Floer homology,
Theorem~\ref{thm:KnotFloerHomology} is still somewhat impractical, as
the chain complex ${C}(\Gamma)$ typically has far too many generators:
for a knot with arc index $n$, the procedure gives a chain complex
with $n!$ generators.  It remains a very interesting challenge to come
up with more efficient methods for calculating the homology of the
complexes we describe here. 

In a different direction, the relationship between our combinatorial
description and Legendrian knots seems tantalizing: one wonders
whether this is perhaps the hint of a connection with the holomorphic
invariants of those objects, compare~\cite{Chekanov},
\cite{Eliashberg}, \cite{Ng}.

We would like to remind the reader that we have kept the introduction
as elementary as possible. The more general results of
Section~\ref{sec:Conclusion} actually lead to a calculation of link
Floer homology for links in $S^3$.  Also, the extra data about the
``knot filtration'' allows one to calculate the concordance invariant
$\tau$ for knots. It is also the input needed to determine the ranks
of Heegaard Floer homology groups of Dehn surgeries on a given knot
$K$, see~\cite{RatSurg}.

\vskip.3cm 
\noindent\textbf{Acknowledgements.} This paper grew out of 
attempts at understanding an earlier preprint by the third author, who
made the revolutionary observation that for Heegaard diagrams of a
certain special form, the corresponding Heegaard Floer homology groups
can be calculated combinatorially. In a different direction, that
preprint also lead to the paper~\cite{SarkarWang}, which gives a
method for describing $\HFa$ of an arbitrary three-manifold in
combinatorial terms.

We are grateful to Matthew Hedden, Mikhail Khovanov, John Morgan, and
Lev Rozansky for their suggestions on an early version of our results.
We are especially grateful to Dylan Thurston for his many interesting
comments, especially for his suggestions for simplifying the Alexander
gradings.  Finally, we owe a great debt of gratitude to Zolt{\'a}n
Szab{\'o}, whose ideas have, of course, had a significant impact on
this present work.  \medskip

\newcommand\D{\mathcal D}
\newcommand\alphas{\mbox{\boldmath$\alpha$}}
\newcommand\betas{\mbox{\boldmath$\beta$}}
\newcommand\ws{\mathbf w}
\newcommand\zs{\mathbf z}

\section {Link Floer homology with multiple basepoints}
\label{sec:Basepoints}

We review here the construction of knot and link Floer homology,
considering the case where the link meets the Heegaard surface in
extra intersection points.  The fact that Heegaard Floer homology can
be extracted from this picture follows essentially from~\cite{Links}.

Let $(\Sigma,\alphas,\betas,\ws,\zs)$ be a Heegaard diagram, where
$\Sigma$ is a surface of genus $g$, $k$ is some positive integer,
$\alphas=\{\alpha_1,...,\alpha_{g+k-1}\}$ are pairwise disjoint,
embedded curves in $\Sigma$ which span a half-dimensional subspace of
$H_1(\Sigma;\Z)$ (and hence specify a handlebody $U_\alpha$ with
boundary equal to $\Sigma$), $\betas=\{\beta_1,...,\beta_{g+k-1}\}$
is another collection of attaching circles specifying $U_\beta$,
and
$\ws=\{w_1,...,w_k\}$ and
$\zs=\{z_1,...,z_k\}$ are distinct marked points with 
\[
\ws,\zs\subset \Sigma-\alpha_1-...-\alpha_{g+k-1}-\beta_1-...-\beta_{g+k-1}.
\]
The data $(\Sigma,\alphas,\betas)$ specifies a Heegaard splitting for
some oriented three-manifold $Y$. In the present applications, we will
be interested in the case where the ambient three-manifold is the
three-sphere, and hence, we make this assumption hereafter.

Let $\{A_i\}_{i=1}^k$ resp. $\{B_i\}_{i=1}^k$ be the connected
components of $\Sigma-\alpha_1-...-\alpha_{g+k-1}$ resp.
$\Sigma-\beta_1-...-\beta_{g+k-1}$.

We suppose that our basepoints are placed in such a manner that each
component $A_i$ or $B_i$ contains exactly one basepoint amongst the
$\ws$ and exactly one basepoint amongst the $\zs$. We can label our
basepoints so that $A_i$ contains $z_i$ and $w_i$, and then $B_i$
contains $w_i$ and $z_{\nu(i)}$, for some permutation $\nu$ of
$\{1,...,k\}$.

In this case, the basepoints uniquely specify an oriented link $L$ in
$S^3=U_\alpha\cup U_\beta$, by the following conventions. For each
$i=1,...,k$, let $\xi_i$ denote an arc in $A_i$ from $z_i$ to $w_i$
and let $\eta_i$ denote an arc in $B_i$ from $w_i$ to
$z_{\nu(i)}$.  Let ${\widetilde\xi}_i\subset U_\alpha$ be an arc
obtained by pushing the interior of $\xi_i$ into $U_\alpha$, and
${\widetilde\eta}_i$ be the arc obtained by pushing the interior of
$\eta_i$ into $U_\beta$.  Now, we can let $L$ be the oriented link
obtained as the sum
$$\bigcup_{i=1}^k\left({\widetilde\xi}_i+{\widetilde\eta}_i\right).$$

\begin{definition}
\label{def:Compat}
In the above case, we say that $(\Sigma, \alphas, \betas, \ws, \zs)$
is a $2k$-pointed Heegaard diagram compatible with the oriented link
$L$ in $S^3$.
\end{definition}

Let $\ell$ denote the number of components of $L$. Clearly, $k\geq
\ell$. In the case where $k=\ell$, these are the Heegaard diagrams
used in the definition of link Floer homology~\cite{Links}, see
also \cite{Knots}, \cite{RasmussenThesis}. In the case where $k>\ell$,
these Heegaard diagrams can still be used to calculate link Floer
homology, in a suitable sense.

\begin{definition}
A {\em periodic domain} is a two-chain of the form
$$P=\sum_{i=1}^k (a_i \cdot A_i + b_i \cdot B_i)$$
which has zero local multiplicity at all of the $\{w_i\}_{i=1}^k$.
A Heegaard diagram is said to be {\em admissible} if every 
non-trivial periodic domain has some positive 
local multiplicities and some negative local multiplicities.
\end{definition}

Consider first the case where our link is in fact a knot.  In this
case, admissibility is automatically satisfied.  Specifically, if we
introduce cyclic orderings of $\{A_i\}_{i=1}^k$ and $\{B_i\}_{i=1}^k$,
$\{w_i\}_{i=1}^k$ and $\{z_i\}_{i=1}^k$, so that $w_i, z_i\in A_i$ and
$w_i, z_{i+1}\in B_i$, then $n_{w_i}(P)=a_i+b_i$ and
$n_{z_i}(P)=a_{i}+b_{i-1}$. The condition that $P$ is a periodic domain
ensures that for each $i$, $a_i+b_i=0$. Thus, if for some $i$ we have
$n_{z_i}(P)>0$ (i.e. $a_{i}+b_{i-1}>0$) then for some other $j$,
$n_{z_j}(P)=a_{j}+b_{j-1}<0$. Conversely, if $n_{w_i}(P)=n_{z_i}(P)=0$
for all $i$, then there is some constant $c$ with all $a_i=c=-b_i$; it
follows readily that $P=0$.

Let $(\Sigma,\alphas,\betas,\ws,\zs)$ be a Heegaard diagram compatible with 
an oriented knot $K$. We will consider Floer homology in the $g+k-1$-fold
symmetric product of the surface $\Sigma$, relative to the pair of 
totally real submanifolds 
\begin{eqnarray*}
\Ta=\alpha_1\times...\times\alpha_{g+k-1}
&{\text{and}}&
\Tb=\beta_1\times...\times\beta_{g+k-1}.
\end{eqnarray*}
Given $\x,\y\in\Ta\cap\Tb$, let $\pi_2(\x,\y)$ denote the space of
homology classes of Whitney disks from $\x$ to $\y$, i.e. maps of the
standard complex disk into $\Sym^{g+k-1}(\Sigma)$ which carry $i$
resp. $-i$ to $\x$ resp. $\y$, and points on the circle with negative
resp. positive real part to $\Ta$ resp. $\Tb$.  (Note that when
$g+k>3$, homology classes of Whitney disks agree with homotopy
classes.)

We consider now the chain complex
$\CFKm(\Sigma,\alphas,\betas,\ws,\zs)$ over the polynomial algebra
$\Field[U_1,...,U_k]$ which is freely generated by intersection points
between the tori ${\mathbb T}_\alpha=\alpha_1\times...\times
\alpha_{g+k-1}$ and ${\mathbb T}_\beta=\beta_1\times...\times\beta_{g+k-1}$ in
${\mathrm{Sym}}^{g+k-1}(\Sigma)$. This module is endowed with the
differential
\begin{equation}
\label{eq:DefD}
\partial^- \x = \sum_{\y\in{\mathbb
    T}_{\alpha}\cap{\mathbb T}_{\beta}} \sum_{\{\phi\in\pi_2({\mathbf
    x},\y)\big| 
\mu(\phi)=1 \}}
\#\left(\frac{{\mathcal
      M}(\phi)}{\mathbb R}\right) U_1^{n_{w_1}(\phi)}\cdot...\cdot
U_k^{n_{w_k}(\phi)}\cdot \y,
\end{equation}
where, as usual, $\pi_2(\x,\y)$ denotes the space of
homology classes of Whitney disks connecting $\x$ to
$\y$, ${\mathcal M}(\phi)$ denotes the moduli space of
pseudo-holomorphic representatives of $\phi$, $\mu(\phi)$ denotes its
formal dimension (Maslov index), $n_{p}(\phi)$ denotes the local
multiplicity of $\phi$ at the reference point $p$ (i.e. the algebraic
intersection number of $\phi$ with the subvariety
$\{p\}\times{\mathrm{Sym}}^{g+k-2}(\Sigma)$), and $\#()$ denotes a
count modulo two.  As usual, in the definition of pseudo-holomorphic
disks, one uses a suitable perturbation of the condition on the disk
that it be holomorphic with respect to the complex structure on
$\Sym^{g+k-1}(\Sigma)$ induced from some complex structure on
$\Sigma$, as explained
in~\cite[Section~\ref{HolDisk:sec:Analysis}]{HolDisk}; see
also~\cite{FloerUnregularized}, \cite{OhFloer},
\cite{FloerHoferSalamon}, \cite{FOOO} for more general discussions. We use here
a sufficiently small perturbation to
retain the property that if $u$ is pseudo-holomorphic, then for all
$p\in\Sigma-\alpha_1-...-\alpha_{g+k-1}-\beta_1-...-\beta_{g+k-1}$,
$n_p(\phi)\geq 0$,
cf.~\cite[Lemma~\ref{HolDisk:lemma:NonNegativity}]{HolDisk}.
For the case where the Heegaard diagram is admissible, it is easy to
see that Equation~\eqref{eq:DefD} gives a finite sum,
compare~\cite{HolDisk}.

The {\em relative Alexander grading} of two intersection points
$\x$ and $\y$ is defined by the formula 
\begin{equation}
  \label{eq:DefRelAlexGrading}
  A(\x)-A(\y)=\left(\sum_{i=1}^n n_{z_i}(\phi)\right)-
  \left(\sum_{i=1}^n n_{w_i}(\phi)\right),
\end{equation}
where $\phi\in\pi_2(\x,\y)$ is any homotopy class from $\x$ to $\y$.
We find it convenient to remove the additive indeterminacy in $A$:
there is a unique choice with the property that
\begin{equation}
  \label{eq:NormalizeRelAlexGrading}
  \sum_{\x\in\Ta\cap\Tb} T^{A(\x)}
  \equiv \Delta_K(T)\cm (1-T^{-1})^{n-1} \ 
  \pmod{2},
\end{equation}
where $\Delta_K(T)$ is the symmetrized Alexander polynomial
which could be made to work over $\Z$ by introducing signs.
(These conventions are chosen to be consistent
with those made in Proposition~\ref{prop:ExtraBasepoints} below.)

Moreover, there is a {\em relative Maslov grading}, defined by
\begin{equation}
        \label{eq:DefRelMasGrading}
M(\x)-M(\y)=\mu(\phi)-2\sum_{i=1}^n n_{w_i}(\phi).
\end{equation}
The relative Maslov grading can be lifted to an absolute grading using
the observation that $(\Sigma,\alphas,\betas,\ws)$ is a
multiply-pointed Heegaard diagram for $S^3$ (a {\em balanced
$n$-pointed Heegaard diagram} in the terminology of~\cite{Links}), and
consequently, if one sets all the $U_i= 0$, the
homology groups of the resulting complex, one obtains a relatively
graded group which is isomorphic to $H_*(T^{k-1};\Field)$
(compare~\cite[Theorem~\ref{Links:thm:InvarianceHFaa}]{Links}).  The
Maslov grading is fixed by the requirement that 
\begin {equation}
        \label{eq:DefAbsMasGrading}
H_*(\CFKm/\{U_i= 0\})\cong H_{*+k-1}(T^{k-1};\Field).
\end {equation}

So far, we have made no reference to the basepoints $\zs$, and indeed,
the complex $\CFm(\Sigma,$ $\alphas,$ $\betas,\ws,\zs)$ so far is the 
chain complex for $\HFm(S^3)$ for a multi-pointed Heegaard diagram, in the
sense of~\cite[Section~\ref{Links:thm:InvarianceHFaa}]{Links}.

This complex admits an {\em Alexander filtration} defined by the
convention that any element $\x\in\Ta\cap\Tb$ has
Alexander filtration level $A(\x)$, and multiplication by the
variables $U_i$ drops Alexander filtration by one, i.e.
$$A(U_1^{a_1}\cdot ...\cdot U_k ^{a_k}\cm \x) =A(\x)-a_1-...-a_k.$$
Non-negativity of local multiplicities of pseudo-holomorphic disks
ensures that this function indeed defines a filtration on the complex;
i.e. we have an increasing sequence of subcomplexes
$\Filt^-(K,m)\subset \CFKm(K)$ indexed by integers $m$, which are
generated over $\Field[U_1,...,U_k]$ by intersection points $\x$ with
$A(\x)\leq m$.

In the case where $k=1$, the above construction gives the chain homotopy
type of the ``knot filtration'' on $\CFm(S^3)$, called $\CFKm(K)$,
which is a chain complex over the polynomial algebra $\Field[U]$.
In this case, it was shown in~\cite{Knots} and~\cite{RasmussenThesis}
that the filtered chain homotopy type of the complex is a knot invariant. 
Our goal here is to show that this filtered chain homotopy type
is also independent of $k$.

In practice, it is often more convenient to consider 
the simpler complex
$C(\Sigma,\alphas,\betas,\ws,\zs)$
generated by intersection points of $\Ta$ and $\Tb$ with coefficients
in $\Field$, endowed with the differential
$$\partial \x = \sum_{\y\in{\mathbb
    T}_{\alpha}\cap{\mathbb T}_{\beta}} \sum_{\{\phi\in\pi_2({\mathbf
    x},\y)\big|
\begin {tiny} 
\begin{array}{c}
\mu(\phi)=1, \\
n_{w_i}(\phi)=n_{z_i}(\phi)=0 ~~~~~~~~ \forall i=1,...,n
\end{array} \end {tiny}\}}
\#\left(\frac{{\mathcal M}(\phi)}{\mathbb R}\right)\cdot \y.  $$ 
For this complex, the function $A$ defines an Alexander grading which
is preserved by the differential. One can think of
$C(\Sigma,\alphas,\betas,\ws,\zs)$ as obtained from
$\CFKm(\Sigma,\alphas,\betas,\ws,\zs)$ by first setting all the
$U_i= 0$, and then taking the graded object associated to the
Alexander filtration. 

The following proposition shows how to extract the usual knot Floer
homology from the above variants using multiple basepoints.  The
result is an adaptation of the results from
\cite[Section~\ref{Links:subsec:SimpleStabilizations}]{Links}, but we
sketch the proof here for the reader's convenience.

\begin{proposition}
  \label{prop:ExtraBasepoints} Let $(\Sigma,\alphas,\betas,\ws,\zs)$
  be a $2k$-pointed Heegaard diagram compatible with a knot
  $K$.  Then, the filtered chain homotopy type of
  $\CFKm(\Sigma,\alphas,\betas,\ws,\zs)$, thought of as a complex over
  $\Field[U]$ where $U$ can be any $U_i$, agrees with the filtered
  chain homotopy type of $\CFKm(K)$.
  Moreover, we have an
  identification \begin{equation} \label{eq:WithoutZs}
  H_*(C(\Sigma,\alphas,\betas,\ws,\zs),\partial)\cong
  {\widehat{HFK}}(K)\otimes V^{\otimes (k-1)}, \end{equation} where
  $V$ is the two-dimensional vector space spanned by two generators,
  one in bigrading $(-1,-1)$, another in bigrading $(0,0)$.
\end{proposition}

We first establish the following:

\begin {lemma}
\label {lemma:Special}
Let $k$ be an integer greater than one.
After a series of isotopies, handleslides and stabilizations, any
$2k$-pointed Heegaard diagram $(\Sigma,\alphas,\betas,\ws,\zs)$
compatible with a knot $K$, can be transformed into one
with the following properties:
\begin{itemize}
\item there are curves $\alpha_1\in\alphas$ and $\beta_1\in\betas$
which bound disks $A_1$ and $B_1$ in $\Sigma$
\item $A_1\cap B_1$ contains the basepoint $w_1$
\item $\alpha_1$ and $\beta_1$ meet transversally in a pair of 
points
\item $\alpha_1$ is disjoint from all $\beta_j$ with $j\neq 1$, and
$\beta_1$ is disjoint from all $\alpha_j$ with $j\neq 1$.
\end{itemize}
\end{lemma}

\begin{figure}
\begin{center}
\mbox{\vbox{\epsfbox{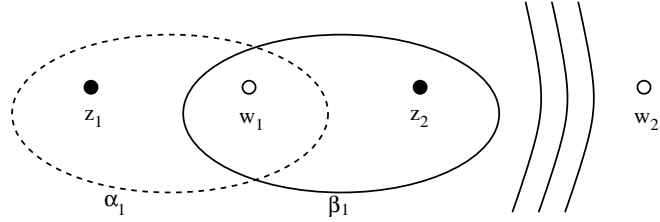}}}
\end{center}
\caption {{\bf Local picture near $w_1$.} We denote
  $\alpha_i$ by dashed and $\beta_j$ by solid lines. The basepoint
  $w_2$ can be connected to $z_2$ by an arc which crosses $\beta_1$,
  and possibly a collection of other $\beta$-circles (but no
  $\alpha$-circles).}
\label{fig:Stabilization}
\end{figure}

\begin {proof} Start from a $2k$ pointed Heegaard diagram
  $(\Sigma,\alphas,\betas,\ws,\zs)$ compatible with $K$.  Let $A_1$ be
  the component of $\Sigma-\alphas_1-...-\alphas_{g+k-1}$ which
  contains $w_1\in\ws$, and $B_1$ be the component of
  $\Sigma-\betas_1-...-\betas_{g+k-1}$ containing $w_1$. In particular
  if $z_1, z_2\in\zs$ are contained inside $A_1$ and $B_1$
  respectively, then $z_1\neq z_2$ (since otherwise the basepoints
  $w_1$ and $z_1$ would determine a closed component of $K$; and since
  $K$ is a knot, we could conclude that $k=1$).  After a sequence of
  handleslides amongst the $\alphas$ and $\betas$ which do not cross
  any of the basepoints $\ws$, $\zs$, we can reduce to the case where
  $A_1$ and $B_1$ are both disks. Let $\alpha_1$ and $\beta_1$ denote
  the boundaries of $A_1$ and $B_1$ respectively.

Note that, since $z_1= A_1\cap \zs$ and $z_2= B_1\cap \zs$ and
$z_1\neq z_2$, the intersection $A_1 \cap B_1$ does not contain any
$z_i\in\zs$. In fact, the various arcs $A_1\cap \beta_1$ divide $A_1$
into a collection of planar regions, one of which contains $w_1$,
another of which contains $z_1$. All the other regions have no
basepoints in them, and we call these {\em unmarked regions}. We can
perform finger moves to eliminate all of the unmarked bigons in $A_1$,
by which we mean unmarked regions in $A_1-A_1\cap \beta_1$ whose
closure meets $\beta_1$ in a single component in $A_1$. Note that this
might involve cancelling also intersection points betweeen $\alpha_1$
with $\beta_j$ for some $j\neq 1$. (See Figure~\ref{fig:FingerMoves}.)
After doing this, $A_1-A_1\cap \beta_1$ consists of one region which
is a bigon marked with $w_1$, another which is a bigon marked with
$z_1$, and some unmarked regions which are all rectangles (i.e. 
their boundary meets $\beta_1$ in two components).

\begin{figure}
\begin{center}
\mbox{\vbox{\epsfbox{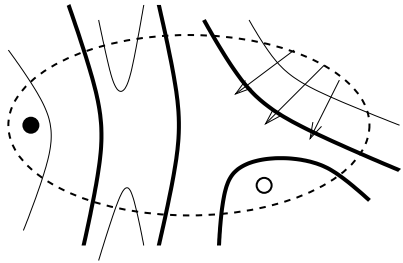}}}
\end{center}
\caption {{\bf Finger moves.}
The circle $\alpha_1$ is indicated by the dashed line. It bounds the
disk $A_1$, which contains the basepoint $w_1$ (indicated by the
hollow dot) and the basepoint $z_1$ (indicated by the dark dot). Other
arcs belong to various $\beta$-circles, which divide $A_1$ into planar
regions, with $\beta_1$ arcs denoted by the thicker lines and other
$\beta_j$ (with $j\neq 1$) by thinner ones.  Performing the finger
move on $\alpha_1$ as indicated by the arrows, we can reduce the
number of unmarked bigon regions in $A_1-A_1\cap \beta_1$.}
\label{fig:FingerMoves}
\end{figure}

Now, $A_1\cap B_1$ consists of a bigon marked with $w_1$ and
a (possibly empty) collection of unmarked rectangles. We can reduce
the number of unmarked rectangular regions in $A_1\cap B_1$ by a
stabilization, followed by four handleslides, as illustrated in
Figure~\ref{fig:ReduceRectangles}.

\begin{figure}
\begin{center}
\mbox{\vbox{\epsfbox{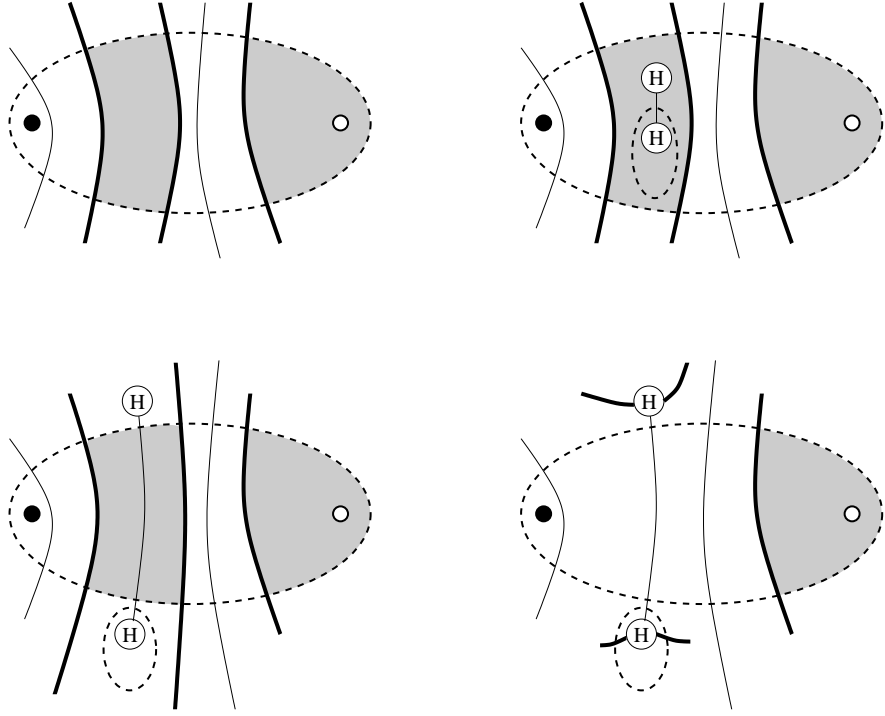}}}
\end{center}
\caption {{\bf Reducing rectangles in $A_1\cap B_1$.}
  The dashed circle represents $\alpha_1$, thick lines represent arcs
  from $\beta_1$, the thin arcs represent arcs from the other
  $\beta_j$ (with $j\neq 1$), and the shaded regions represent
  $A_1\cap B_1$. We eliminate the rectangular region in $A_1\cap B_1$
  by first stabilizing as in the second picture, introducing a new
  handle, represented by the two circles marked with $H$, along with
  the new dashed $\alpha$-circle $\alpha_2$, and the new $\beta$-circle
  $\beta_2$ indicated in the second picture by the thin arc running
  through the handle. Handlesliding $\alpha_1$ over $\alpha_2$ twice,
  we obtain the third picture.  Handlesliding $\beta_1$ over $\beta_2$
  twice, we end up with the fourth picture, which has one fewer
  (rectangular) component in $A_1\cap B_1$.}
\label{fig:ReduceRectangles}
\end{figure}

Finally, by performing handleslides of the additional
$\alpha_i$ and $\beta_j$ over $\alpha_1$ and $\beta_1$ respectively
(followed by some isotopies), we can arrange for the circles
$\alpha_1$ and $\beta_1$ to be disjoint from all the other $\alpha_i$
and $\beta_j$. 
\end {proof}

\vskip.3cm
\noindent{\bf{Proof of Proposition~\ref{prop:ExtraBasepoints}.}}
We use induction on $k$. Obviously, in the case where $k=1$, there is nothing to prove.

When $k > 1,$ using Lemma~\ref{lemma:Special} we can reduce to the case where our diagram 
$$(\Sigma,\alphas,\betas,\ws,\zs)=(\Sigma,\{\alpha_1,...,\alpha_{g+k-1}\},\{\beta_1,...\beta_{g+k-1}\}, 
\{w_1,...,w_k\}, \{z_1,...,z_k\}),$$
has the special form in Figure~\ref{fig:Stabilization}. Note that under all the Heegaard moves used in Lemmaa~\ref{lemma:Special}, the filtered chain
homotopy type of the associated chain complex remains invariant, as
in~\cite{HolDisk}, \cite{Knots}, \cite{RasmussenThesis}.

We can de-stabilize the original diagram to get a $k-1$-pointed Heegaard
diagram 
$$(\Sigma,\{\alpha_2,...,\alpha_{g+k-1}\},
\{\beta_2,...,\beta_{g+k-1}\},\{w_2,...,w_{k}\},\{z_2,...,z_k\})$$
for the same knot.
Let $C'$ denote its corresponding filtered Heegaard Floer complex
$$\CFKm(\Sigma,\{\alpha_2,...,\alpha_{g+k-1}\},
\{\beta_1,...,\beta_{g+k-1}\},\{w_2,...,w_{k}\},\{z_2,...,z_k\}),$$
thought of as a module over the polynomial algebra
$\Field[U_2,...,U_k]$.  It is generated by the corresponding
intersection points $X'$ of the tori $\Ta'$ and $\Tb'$ in
$\Sym^{g+k-2}(\Sigma)$.

  The set of generators $X$ 
  of the complex $\CFKm(\Sigma,\alphas,\betas,\ws,\zs)$ 
  has the form $X'\times\{x,y\}$, where $x$ and $y$ are the two 
  points of intersection of $\alpha_1$ and $\beta_1$. Let $C_x$
  be the subgroup of $C$ generated by intersection
  points of type $X'\times \{x\}$ and $C_y$ be the
  subgroup generated by those of type $X'\times\{y\}$. 
  
  It is shown in the proof
  of~\cite[Proposition~\ref{Links:prop:StabilizationInvarianceOne}]{Links}
  that for a suitable choice of complex structure on $\Sigma$, the
  chain complex $\CFKm(\Sigma,\alphas,\ws,\zs)$ is identified with the
  mapping cone of the chain map 
  \begin{equation}
    \label{eq:MappingCone}
    U_1-U_2\colon C'[U_1]  \longrightarrow C'[U_1],
  \end{equation}
  where the domain is identified with $C_x$
  and the range with $C_y$. Specifically,  under the natural
  identifications of groups $C_x\cong C'[U_1]$, $C_y\cong C'[U_1]$, we
  have that the differential of $C$ is identified with the matrix
  $$\left(
    \begin{array}{ll}
      \partial' & {U_1-U_2} \\
      0 & \partial'
    \end{array}
  \right),
  $$
  where here the variable $U_2$ corresponds to the basepoint $w_2$
  which lies in the same
  $\Sigma-\alpha_1-...-\alpha_{g+k-1}$-component as $z_2$.  In this
  mapping cone, the Alexander filtration of a generator in $C'[U_1]$
  thought of as supported in $C_x$ is one higher than the Alexander
  filtration of the corresponding element, thought of as supported in
  $C_y$.
  
  Let us say a few more words about the identification of
  $\CFKm(\Sigma,\alphas,\betas,\ws,\zs)$ with the mapping cone used
  above, referring the interested reader
  to~\cite[Proposition~\ref{Links:prop:StabilizationInvarianceOne}]{Links}
  for more details. Think of $\Sigma$ as formed by the connected sum
  of $\Sigma'$ with a genus zero surface $S$ containing both
  $\alpha_1$ and $\beta_1$.  Fixing conformal structures on $\Sigma'$
  and the sphere, we obtain a one-parameter family of conformal structures on
  $\Sigma$ by inserting a connected sum tube isometric to $[0,T]\times
  S^1$, and allowing $T$ to vary. When $T$ is sufficiently large, the
  chain complex $\CFKm(\Sigma,\alphas,\ws,\zs)$ can be identified with
  the mapping cone of
  $$U_1-f\colon C'[U_1]\longrightarrow C'[U_1],$$
  where $f$ is a map which counts points in a fibered product of
  moduli spaces of disks coming from $\Sigma$ and $S$, fibered over a
  non-trivial symmetric product of the disk, where the maps are
  obtained as the preimage of the connected sum points $p$ and $q$ in
  $\Sigma'$ and $S$.  (The term $U_1$ fits into this picture formally
  as the fibered product over the empty symmetric product.)  To
  understand $f$ (and identify it with $U_2$), we must consider a
  second parameter $s$ in the space of conformal structures on
  $\Sigma$, which is given by moving the connected sum point $q$ in
  $S$. Indeed, it will be useful to move the connected sum point $q\in
  S$ towards $\alpha_{1}$ (so that as $s\mapsto \infty$, $q$ limits
  onto $\alpha_1$). In fact, for $q$ sufficiently close to $\alpha_1$
  (i.e. $s$ sufficiently large), the only non-empty moduli space which
  contributes to this fiber product consists of holomorphic disks in
  $\Sym^{g+k-2}(\Sigma)$ with Maslov index equal to two which carry
  some fixed point $m$ in the disk (whose distance to the
  $\alpha$-boundary of the disk goes to zero as $s\mapsto \infty$)
  into $q\times \Sym^{g+k-3}$. When $m$ is sufficiently close to the
  $\alpha$-boundary, the count of these disks is identified with the
  count of Maslov index two $\alpha$-boundary degenerations for
  $\Sigma$ with local multiplicity $1$ at the connected sum point $p$.
  Of course, the contribution of these boundary degenerations is given
  by multiplication by $U_2$. This gives the identification of
  $\CFKm(\Sigma,\alphas,\betas,\ws,\zs)$ with the mapping cone of
  Equation~\eqref{eq:MappingCone}.  (Note that we broke the symmetry
  in the construction by moving $q$ towards $\alpha_1$ rather than
  $\beta_1$. If we moved $q$ towards $\beta_1$ instead, we would
  identify $\CFKm(\Sigma,\alphas,\betas,\ws,\zs)$ with the mapping cone
  of $U_1-U_{k}$.)

  For the second assertion of the proposition, view all the $U_i$ as
  being set to zero, and the above argument shows that
  $$C(\Sigma,\alphas,\betas,\ws,\zs)\cong
  C(\Sigma',\alphas',\betas',\ws',\zs')\otimes V;$$ and hence a
  corresponding identification holds on the level of homology.
  Iterating this until we remain with two basepoints, we obtain the
  stated identification. \qed
\medskip

It is perhaps more traditional to consider the filtration of
$\CFa(S^3)$ (rather than $\CFm$). This filtration is induced from the
filtration of $\CFm(S^3)$ by setting $U=0$.  According to
Proposition~\ref{prop:ExtraBasepoints}, this filtration is obtained
as the induced filtration of 
$\CFKm(\Sigma,\alphas,\betas,\ws,\zs)/U_1$.

\subsection{Modifications for links} 
\label{subsec:Links}

Recall that knot Floer homology has a generalization to the case of
oriented links $\orL$. For an $\ell$-component, oriented link
$\orL$ in the three-sphere, this
takes the form of a multi-graded theory
$$\HFLa(\orL)=\bigoplus_{d\in\Z,h\in\H}\HFLa_d(\orL,h),$$
where $\H\cong H_1(S^3-\orL)\cong \Z^\ell$,
with the latter isomorphism induced by an ordering of the link components.
We sketch now the changes to be made to the above discussion 
for the purposes of understanding link 
Floer homology for Heegaard diagrams with extra basepoints,
where by ``extra'' here we mean more than twice $\ell$.

Suppose now that $(\Sigma,\alphas,\betas,\ws,\zs)$ is a Heegaard
diagram compatible with an oriented link $\orL$ in the sense of
Definition~\ref{def:Compat}.

We find it convenient to label the basepoints keeping track of which
link component they belong to. Specifically, suppose $L$ is a link with $\ell$
components, and for $i=1,...,\ell$,
we choose $k_i$ basepoints to lie on the $i^{th}$ component.
Letting $S$ be the index set of pairs $(i,j)$ with $i=1,...,\ell$ and 
$j=1,...,k_i$. We now have basepoints
$\{z_{i,j}\}_{(i,j)\in S}$ and $\{w_{i,j}\}_{(i,j)\in S}$.

We can now form the chain complex
$\CFLm(\Sigma,\alphas,\betas,\ws,\zs)$ defined over
$\Field[\{U_{i,j}\}_{(i,j)\in S}]$ analogous to the version before,
generated by intersection points of $\Ta\cap\Tb$, with differential
\[\partial^- \x = \sum_{\y\in{\mathbb
    T}_{\alpha}\cap{\mathbb T}_{\beta}} \sum_{\{\phi\in\pi_2({\mathbf
    x},\y)\big| \mu(\phi)=1\}}
\#\left(\frac{{\mathcal
      M}(\phi)}{\mathbb R}\right)\cm \left(\prod_{(i,j)\in S} U_{i,j}^{n_{w_{i,j}}(\phi)}\right) \cm \y.
\]

This complex has a relative Maslov grading as before. It also
has a relative Alexander grading which in this case is an $\ell$-tuple
of integers, 
$$A\colon \Ta\cap\Tb \longrightarrow \Z^{\ell},$$
determined up to an overall additive constant by the formula
$$A({\mathbf
x})-A(\y)=\left(\sum_{j=1}^{k_1} (n_{z_{1,j}}(\phi)-
n_{w_{1,j}}(\phi)), 
...,\sum_{j=1}^{k_\ell} (n_{z_{\ell,j}}(\phi)-n_{w_{\ell,j}}(\phi))
\right).
$$
The indeterminacy in this case is a little more unpleasant to pin
down (i.e. one must go beyond the multi-variable Alexander polynomial,
which can be identically zero), but one can do this with the help of
Proposition~\ref{prop:ExtraBasepointsLink}.

The complex $\CFLm(\Sigma,\alphas,\betas,\ws,\zs)$ inherits an
Alexander filtration induced by the Alexander multi-grading of
$\Ta\cap\Tb$, and the convention that $U_{i,j}$ drops the
multi-grading by the $i^{th}$ basis vector.
In the case where $k=\ell$, the filtered chain homotopy type
of $\CFLm(\Sigma,\alphas,\betas,\ws,\zs)$ was shown to be
a link invariant in~\cite{Links}; it is the link filtration
$\CFLm(\orL)$.

The analogue of $C(\Sigma,\alphas,\betas,\ws,\zs)$ can be defined 
as well: it is generated by intersection points of $\Ta$ and
$\Tb$ over $\Field$, endowed with the differential
$$\partial \x = \sum_{\y\in{\mathbb
    T}_{\alpha}\cap{\mathbb T}_{\beta}} \sum_{\{\phi\in\pi_2({\mathbf
    x},\y)\big| 
\begin {tiny}
\begin{array}{c}
\mu(\phi)=1, \\
n_{w_{i,j}}(\phi)=n_{z_{i,j}}(\phi)=0 ~~~~~~~~ \forall (i,j)\in S
\end{array}
\end {tiny}
\}}
\#\left(\frac{{\mathcal M}(\phi)}{\mathbb R}\right)\cdot \y.  $$ 
This differential drops Maslov grading by one and
preserves the Alexander multi-grading,
and hence the homology groups $H_*(C(\Sigma,\alphas,\ws,\zs))$
inherit a Maslov grading and an Alexander multi-grading.

\begin{proposition}
  \label{prop:ExtraBasepointsLink} Let
  $(\Sigma,\alphas,\betas,\ws,\zs)$ be a $2k$-pointed admissible
  Heegaard diagram compatible with an oriented link $\orL$, with $k_i$
  pairs of basepoints corresponding to the $i^{th}$ component of
  $\orL$.  Then, there is a filtered chain homotopy equivalence
  $\CFLm(\Sigma,\alphas,\betas,\ws,\zs)$ with the usual link
  filtration $\CFLm(\orL)$, viewed as a chain complex over
  $\Field[\{U_{i,j}\}_{(i,j)\in S}]$.   Moreover, there are
  (relative) multi-graded identifications \[
  H(\Sigma,\alphas,\betas,\ws,\zs)\cong
  {\widehat{HFL}}(\orL)\otimes \bigotimes_{i=1}^\ell V_i^{\otimes
    (k_i-1) }, \] where $V_i$ is the two-dimensional vector space
  spanned by one generator in Maslov and Alexander gradings zero, and
  another in Maslov grading $-1$ and Alexander grading corresponding to
  minus the $i^{th}$ basis vector. 
\end{proposition}

\begin{proof}
  This follows as in the proof of
  Proposition~\ref{prop:ExtraBasepoints}, with a little extra care
  taken to ensure that all Heegaard diagrams remain admissible while 
performing Heegaard moves, as in \cite[Proposition 7.2]{HolDisk}.  

 \end{proof}

\section{Proof of Theorem~\ref{thm:KnotFloerHomology} and its generalizations}
\label{sec:Conclusion}

The alert reader will have noticed by now that the toroidal grid
diagrams from the introduction are a special case of the
multiply-pointed Heegaard diagrams from Section~\ref{sec:Basepoints}.
The Heegaard surface is the torus $\Torus$, the $\alpha$-circles are
the horizontal circles, and the $\beta$-circles are the vertical ones.
The basepoints $\{w_i\}_{i=1}^n$ are the white dots, and
$\{z_i\}_{i=1}^n$ are the black ones. Since for each $i$ and $j$,
$\alpha_i$ and $\beta_j$ intersect in the single point $(i,j)$, we see
that the generators $X$ are, of course, the intersection points of
$\Ta$ with $\Tb$ in $\Sym^n(\Torus)$. In our coordinate system, these
generators can be thought of as graphs of permutations on $n$ letters.
To apply the results from Section~\ref{sec:Basepoints}, we must
verify the following:

\begin{lemma}
  \label{lemma:AlexanderGrading}
  The Alexander grading of generators $X$, as specified by the
  Heegaard diagram given by the grid diagram (characterized by the
  Equations~\eqref{eq:DefRelAlexGrading}
  and~\eqref{eq:NormalizeRelAlexGrading}) coincides with the function
  $A\colon X\longrightarrow \Z$ defined in the introduction
  (Equation~\eqref{eq:Alexanders}).
\end{lemma}

\begin{proof}
  Let $A'$ denote the Alexander grading of generators specified by the
  Heegaard diagram, and let $A$ denote the function defined in the
  introduction. Our aim is to show that $A=A'$.
  
  Recall that the Alexander grading $A'$ is determined up to an
  overall additive constant by the formula $$A'({\mathbf
    x})-A'({\mathbf y})=\left(\sum_{i=1}^n n_{z_i}(\phi)\right)-
  \left(\sum_{i=1}^n n_{w_i}(\phi)\right),$$
  where
  $\phi\in\pi_2(\x,\y)$ is any homology class connecting $\x$ to $\y$.
  The right-hand side of this equation can be interpreted as the
  oriented intersection number of the knot $K$ with the two-chain
  associated to $\phi$; or alternatively as the linking number of $K$
  with $\gamma_{\x,\y}=\partial{\mathcal D}(\phi)$ (which we think of
  now as an embedded curve in the three-sphere supported near its
  Heegaard torus).  We can also think of this linking number as the
  intersection number of a Seifert surface for $K$ with
  $\gamma_{\x,\y}$. Deforming $\gamma_{\x,\y}$ (without changing its
  intersection number with the Seifert surface for $K$) so that the
  horizontal segments are far under the Heegaard surface, and the
  vertical ones are far above it (so that each $x_i\in\x$ is the
  projection of an arc in $\gamma_{\x,\y}$ which points vertically
  downwards, while each $y_i\in\y$ is the projection of an arc in
  $\gamma_{\x,\y}$ which points vertically upwards), we can arrange
  that all the intersection points of $\gamma_{\x,\y}$ with the
  Seifert surface occur in the arcs over $x_i$ and $y_i$.  Thus, we
  have established that for any two generators $\x,\y\in X$,
  \begin{equation}
    \label{eq:AlexUpToFactor}
    A(\x)-A(\y)=A'(\x)-A'(\y),
  \end{equation}
  or equivalently, that there is some $\kappa$ with the property that
  for any generator $\x\in X$,
  $A(\x)=A'(\x)+\kappa$.
  
  The proof that $\kappa=0$ is elementary, albeit tedious. We sketch
  it here, leaving the details as an exercise for the interested
  reader; compare also~\cite{MOST}. One first checks that $\kappa$ is a knot invariant, by
  verifying that it is unchanged by vertical and horizontal rotations
  of the toroidal grid diagram, as well as by the Reidemeister moves
  from ~\cite{Dynnikov} (see also~\cite{Cromwell}), which relate any
  two planar grid diagrams of the same knot. Then it suffices to show
  that the rational function of $T$ determined by the expression $$
  Q(K)=\frac{\sum_{\x\in X} T^{A(\x)}}{(1-T^{-1})^{n-1}}, $$ which we
  know is $T^{\kappa}\cm \Delta_K(T)$, is actually a symmetric Laurent
  polynomial in $T$. This can be done, for example, by verifying that
  it agrees with the symmetrized Alexander polynomial modulo two,
  using the skein relation: $$Q(K_+)-Q(K_-)\equiv
  (T^{\OneHalf}-T^{-\OneHalf})\cm Q(K_0)\pmod{2}$$ The skein relation
  for $Q$ can be readily verified by realizing the skein moves in grid
  position.  Finally, a straightforward calculation in a $2\times 2$
  diagram shows that $Q(K)=1$ when $K$ is the unknot.  It follows that
  $Q(K)$ is the symmetrized Alexander polynomial modulo two, and in particular
  that $Q$ is symmetric.

\end{proof}

\begin{lemma}
  \label{lemma:MaslovGrading}
  The Maslov grading of generators $X$, as specified by the Heegaard
  diagram $(\Torus,\{\alpha_1,...,\alpha_n\},\{\beta_1,...,\beta_n\},
  \{w_1,...,w_n\},\{z_1,...,z_n\})$ and characterized by
  Equations \eqref{eq:DefRelMasGrading}
  and \eqref{eq:DefAbsMasGrading}, coincides with the function
  $M\colon X \longrightarrow \Z$ defined in the introduction
  (characterized by Equation~\eqref{eq:RelativeMaslovGrading} and the
  normalization that $M(\x_0)=1-n$).
\end{lemma}

\begin{proof}
 For $\phi \in \pi_2(\x, \y),$ we claim that its Maslov index $\Mas(\phi)$ 
is given by the formula
\begin {equation}
\label {eq:lipshitz}
\Mas(\phi)=P_\x(\cald(\phi))+P_\y(\cald(\phi)).
\end {equation}

This is a particular case of Lipshitz's formula for the Maslov index in an 
arbitrary Heegaard diagram~\cite{LipshitzCyl}. However, for domains on a 
grid diagram, we can also give an elementary proof as follows.

Let $\Mas'(\phi)$ denote the quantity on the right-hand side of Equation 
\eqref{eq:lipshitz}. First, note that $\Mas'(\phi) = 1$ when the domain 
$\cald(\phi)$ associated to a homology class $\phi\in\pi_2(\x,\y)$ is a 
rectangle $r$ in the torus which contains none of the components of $\x$ 
in its interior; in this case we also have $\Mas(\phi) = 1,$ because the 
moduli space of complex structures on a disk with four marked points on 
the boundary is one-dimensional.

Next, consider the natural map given by juxtaposition of flow lines:
$$ * : \pi_2(\x, \y) \times \pi_2(\y, \z) \to \pi_2(\x,\z).$$

The Maslov index is additive under this operation, i.e. $\Mas(\phi_1 * 
\phi_2) = \Mas(\phi_1) + \Mas(\phi_2).$ We claim that the same is true for 
$\Mas'.$ Indeed, the relation
$$ P_\x(\cald(\phi_1)) +P_\y(\cald(\phi_1))+  P_\y(\cald(\phi_2)) 
+P_\z(\cald(\phi_2)) = P_\x(\cald(\phi_1 * \phi_2))+ P_\z(\cald(\phi_1 * 
\phi_2))$$
is equivalent to
\begin {equation} 
\label {eq:linking}
P_\x(\cald(\phi_2)) - P_\y(\cald(\phi_2)) = P_\y(\cald(\phi_1)) - 
P_\z(\cald(\phi_1)).
\end {equation}

Let us denote by $\gamma_{\x, \y}^{NE}, \gamma_{\x, \y}^{NW}, \gamma_{\x, 
\y}^{SW},\gamma_{\x, \y}^{SE}$ small translates of the curve $\gamma_{\x, 
\y}= \partial \cald(\phi_1)$ on the torus, in the four diagonal 
directions. As in the proof of Lemma~\ref{lemma:AlexanderGrading}, we 
deform these curves by pushing their horizontal arcs under the Heegaard 
surface and their vertical arcs above the Heegaard surface; thus we can 
think of them as embedded curves in the three-sphere. The left-hand side 
of Equation~\eqref{eq:linking} is then the average of the intersection 
numbers of the surface $\cald(\phi_2)$ with each of $\gamma_{\x, \y}^{NE}, 
\gamma_{\x, \y}^{NW}, \gamma_{\x, \y}^{SW},\gamma_{\x, \y}^{SE}.$ 
Alternatively, it can be viewed as the average linking number of 
$\gamma_{\y, \z} = \partial \cald(\phi_2)$ with these four curves. 
The right-hand side of \eqref{eq:linking} has a similar 
interpretation. Since the linking number is symmetric, the two sides are 
equal. Therefore, $\Mas'$ is additive under juxtaposition of flow lines.

Now, given an arbitrary pair $\x,\y\in X$, it is easy to construct a
sequence of generators $\x_1,...,\x_m\subset X$ with $\x=\x_1$,
$\y=\x_m$, and $\phi_i\in\pi_2(\x_i,\x_{i+1})$ with the property that
$\cald(\phi_i)$ is a rectangle with no components of $\x_i$ in its
interior. It follows that if we let $\psi=\phi_1*...*\phi_m$, then
$M(\psi)=M'(\psi)$. The alpha curves cut the torus into $n$ annuli
$\{A_i\}_{i=1}^n$, and similarly the beta curves cut it into annuli
$\{B_i\}_{i=1}^n$. The homology classes $\psi$ and $\phi$ differ by
adding or subtracting some number of copies of annuli $A_i$ or $B_j$
(thought of as elements of $\pi_2(\x,\x)$), for which
$\Mas(A_i)=\Mas'(A_i)=2$ (because $A_i$ can be decomposed as a
juxtapositon of two rectangles).  It follows that
$\Mas(\phi)=\Mas(\phi')$.

We have verified Equation \eqref{eq:lipshitz}. It follows now that the 
relative Maslov grading from Equation~\eqref{eq:DefRelMasGrading} 
specializes to Equation~\eqref{eq:RelativeMaslovGrading}.

\begin{figure}
\begin{center}
\mbox{\vbox{\epsfbox{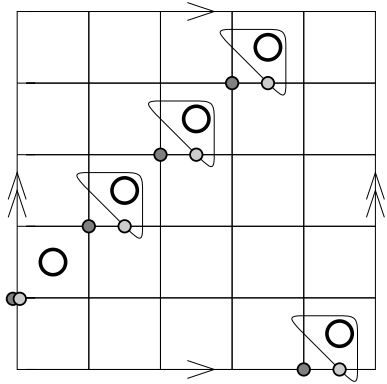}}}
\end{center}   
\caption {{\bf Fixing the Maslov grading.} Handleslide the
vertical circles from left to right, to obtain the smaller null-homotopic
circles $\beta_i'$ encircling the various $w_i$. There is a collection of
triangles connecting the generator $\x_0$, indicated here with the
darkly shaded circles, with the bottom-most generator of $\Ta\cap\Tb'$,
indicated here with the lightly shaded circles.}
\label{fig:VerifyMaslovZero}
\end{figure}
  
  We can lift from the relative to the absolute Maslov grading by
  performing handleslides on the $\beta$-circles in our diagram which
  now are allowed to cross the $z_i$, to reduce to a diagram which has
  $2^{n-1}$ intersection points in $\Ta\cap \Tb'$, and for which all 
  the differentials in the chain complex vanish; indeed, it is       
  identified with the homology of an $n-1$-dimensional torus. The    
  handleslides are performed by successively handlesliding $\beta_i$ 
  over $\beta_{i+1}$ for $i=1,...,n-1$, as pictured in
  Figure~\ref{fig:VerifyMaslovZero}.  It is easy to see now that the 
  generator $\x_0$ from the introduction can be connected to the
  bottom-most generator of the new chain complex by a collection of
  (Maslov index zero) triangles. According to
  (\ref{eq:DefAbsMasGrading}), the grading of $\x_0$ should be $1-n.$
\end{proof}

We can now turn to the following:

\vskip.3cm
\noindent{\bf{Proof of Theorem~\ref{thm:KnotFloerHomology}.}}
The fact that the Alexander and Maslov gradings are identified has
been verified in Lemmas~\ref{lemma:MaslovGrading} and
\ref{lemma:AlexanderGrading} above. It remains to identify the
differentials.

The circles $\alpha_1,...,\alpha_n$ and $\beta_1,...,\beta_n$ cut up
$\Torus$ into $n^2$ squares $D_{i,j}$ with $1\le i,j\le n$. According
to \cite[Proposition~\ref{HolDisk:prop:WhitneyDisks}]{HolDisk},
homology classes of Whitney disks $\phi\in\pi_2(\x,\y)$ are determined
by their underlying two-chain
$$\cald(\phi)=\sum_{i,j}
a_{i j} D_{i,j},$$ where here $a_{i,j}=n_{p_{i,j}}(\phi)$ for some
point $p_{i,j}\in D_{i,j}$. Indeed, if $\x$ and $\y$ correspond to permutations
$\sigma$ and $\tau$, then these induced two-chains are the ones that satisfy
the property for all $i=1,...,n$ that 
$$\partial (\partial {\mathcal D}(\phi)\cap \alpha_i)=(i,\tau(i))-(i,\sigma(i)).$$

To understand the differential, we must count holomorphic 
disks in $\ModFlow(\phi)$ with $\Mas(\phi)=1$. First, we
classify all non-negative homology classes $\phi$ with Maslov index one.

Let $D=\cald(\phi)$.
First observe that if $\partial D$ is $0$ on $(n-1)$ $\alpha$
circles, then it is in fact $0$ on all the $\alpha$ circles, and
$D$ is generated by the annular regions cut out by the $\beta$
circles. Now, if such a thing happens then $\x=\y$, and
its Maslov index is even.

Thus we can assume that $\partial D$ is non-zero on at least two
$\alpha$ circles (say $\alpha_{j_1}$ and $\alpha_{j_2}$) and similarly
non-zero on at least two $\beta$ circles ($\beta_{i_1}$ and
$\beta_{i_2}$).  It follows that there are permutations $\sigma$ and $\tau$
such that
\begin{eqnarray*}
P_{\x}(D)&\geq& p_{(i_1,\sigma(i_1))}(D)+p_{i_2,\sigma(i_2)}(D)\geq 1/2 \\
P_{\y}(D)&\geq& p_{(i_1,\tau(i_1))}(D)+p_{i_2,\tau(i_2)}(D)\geq 1/2.
\end{eqnarray*}
Since $\Mas(\phi)=P_{\x}(D)+P_{\y}(D)=1$, equality must hold
throughout. It follows that
$\partial D$ is non-zero precisely on $\alpha_{j_1},
\alpha_{j_2}, \beta_{i_1}$ and $\beta_{i_2}$, and $D$ is one of the
two rectangles with four vertices $(i_1,j_1)$, $(i_1,j_2)$,
$(i_2,j_1)$, and $(i_2,j_2)$.  Without loss of generality, assume
$\sigma(i_1)= j_1$ and $\sigma(i_2)= j_2$. Then $\tau(i_1)=j_2$ and
$\tau(i_2)=j_1$, and it agrees with $\sigma$ on the rest of the
values. Also for the Maslov index requirement,
$(i,\sigma(i))=(i,\tau(i))$ does not lie in the interior of $D$ for
any other $i$.

Thus, we have established that the only $\phi\in\pi_2(\x,\y)$ with
non-negative local multiplicities and Maslov index equal to one are
those whose underlying domain $r$ is a rectangle of the form $r\in
R_{\x,\y}$ with $P_\x(r)+P_\y(r)=1$.  Moreover, we claim that in this
case, the number of pseudo-holomorphic representatives of $r$ is
odd. In fact, this can be seen by elementary complex analysis, using a
(classical) complex structure on the symmetric product of $\Torus$,
where one shows that in fact the moduli space consists of a single
representative. Indeed, for this choice, the moduli space
$\ModFlow(r)/\R$ can be seen to correspond to involutions of $r$ (with
the complex structure it inherits from $\Torus$) which switch opposite
sides of the rectangle. It is a simple exercise in conformal geometry that
for any rectangle, there is a unique such involution.

We have thus completed the verification that the complex $C(\Gamma)$ from the introduction coincides with the Heegaard Floer
complex $C(\Torus,\alphas,\betas,\ws,\zs)$ in the notation
of Section~\ref{sec:Basepoints}. Theorem~\ref{thm:KnotFloerHomology}
now follows directly from Proposition~\ref{prop:ExtraBasepoints}
(Equation~\eqref{eq:WithoutZs}).
\qed

\subsection{Other variants}
There are other variants of Theorem~\ref{thm:KnotFloerHomology}, which
should be clear from the constructions thus far. We state several of
them for completeness.

Label the white dots $\{w_1,...,w_n\}$, and let $n_{w_i}(r)$
denote the local multiplicity of $r$ at $w_i$.
Consider the chain complex $C^-(\Gamma)$ over the algebra
$\Field[U_1,...,U_n]$ also generated by $X$, endowed with the differential
$${\partial}^- \x = \sum_{\y\in X} \sum_{r\in R_{\x,\y}}
\left\{\begin{array}{ll} 1 & {\text{if $P_{\x}(r)+P_{\y}(r)=1$}} \\ 0
    & {\text{otherwise}}
\end{array}
\right\} U_1^{n_{w_1}(r)}\cm...\cm U_n^{n_{w_n}(r)} \cm \y,$$
thought
of as a filtered chain complex where the filtration level of each
generator $\x\in X$ is its Alexander grading, and multiplication by
the variable $U_i$ drops filtration level by one.

\begin{theorem}
  \label{thm:KnotFloerHomologyMinus} Fix a grid presentation $\Gamma$ of
  a knot $K$, with grid number $n$. The filtered chain homotopy
  type of $K$ coincides with the filtered chain homotopy
  type of the knot filtration $\CFm(S^3,K)$.
\end{theorem}

\begin{proof}
  The proof of Theorem~\ref{thm:KnotFloerHomology} identifies the
  filtered chain complex $C^-(\Gamma)$ with the complex denoted
  $\CFKm(\Torus,\alphas,\betas,\ws,\zs)$ in
  Section~\ref{sec:Basepoints} which, by
  Proposition~\ref{prop:ExtraBasepoints}, is identified with
  $\CFKm(K)$.
\end{proof}

Of course, the other filtrations $\CFKinf(S^3,K)$ and $\CFKp(S^3,K)$
from~\cite{Knots} can be extracted from this information.

We call attention to another other construction, which gives a concordance
invariant $\tau(K)$ for knots~\cite{4BallGenus},
\cite{RasmussenThesis}.  This is a homomorphism from the smooth
concordance group of knots to the integers, which can be used to bound
the four-ball genus of knots, giving an alternate proof of the theorem
of Kronheimer and Mrowka~\cite{KMMilnor} confirming Milnor's
conjecture for the unknotting numbers of torus knots. This feature
underscores its similarity with Rasmussen's concordance invariant
$s(K)$~\cite{RasmussenSlice} from Khovanov homology~\cite{Khovanov}.
However, these two invariants are known to be linearly
independent~\cite{HeddenOrding}.

Recall that the filtration $\CFKm(K)$
of $\CFm(S^3)$ induces also a filtration $\{{\widehat {\mathcal
    F}}_m(S^3)\}_{m\in\Z}$ of $\CFa(S^3)=\CFm(S^3)/U\cdot\CFm(S^3)$. The
concordance invariant $\tau(K)$ is by definition the minimal $m\in\Z$
with the property that the map
$$H_*({\widehat {\mathcal F}}_m(S^3))\longrightarrow \HFa(S^3)\cong \Field
$$
is non-trivial.

We have a corresponding chain complex ${\widehat
  C}(\Gamma)=\CFm(\Gamma)/U_1\cdot\CFm(\Gamma)$; i.e.  whose differential
is given by
$${\widehat \partial} \x = \sum_{\y\in X} \sum_{r\in R_{\x,\y}}
\left\{\begin{array}{ll} 1 & {\text{if $P_{\x}(r)+P_{\y}(r)=1$ and
$n_{w_1}(r)=0$}} \\ 0 & {\text{otherwise}}
\end{array}
\right\} U_2^{n_{w_2}(r)}\cm...\cm U_n^{n_{w_n}(r)} \cm \y,$$

This is equipped with subcomplexes
${\widehat F}(K,m)\subset {\widehat C}(\Gamma)$,
generated by elements $U_2^{a_2}\cm...\cm U_n^{a_n}\cm\x$
with integral $a_i\geq 0$, and 
$$A(\x)-a_2-...-a_n\leq m.$$

\begin{corollary}
  The concordance invariant $\tau(K)$ is the minimal $m$ for which
  the map induced on homology 
  $$i_*\colon H_*({\widehat {\mathcal F}}(K,m))\longrightarrow H_*({\widehat C}(\Gamma))$$
  is non-trivial.
\end{corollary}

\begin{proof}
  Theorem~\ref{thm:KnotFloerHomologyMinus} actually gives an
  identification of the filtered chain homotopy type of $\CFKa(K)$
  with ${\widehat C}(\Gamma)$. The result then follows from the definition
  of $\tau(K)$.
\end{proof}

Consider now the case of link Floer homology. In order to use the
Heegaard diagram associated to a grid diagram to calculate link
Floer homology, we must verify that it is admissible.

\begin {lemma}
  The diagram $(\Torus,\{\alpha_1,...,\alpha_n\},\{\beta_1,...,\beta_n\},
  \{w_1,...,w_n\},\{z_1,...,z_n\})$ is admissible.
\end {lemma}

\begin{proof}
  The formal differences $A_i - B_i$ span the space of periodic
  domains. Drawing $\Torus$ as a square, with equally spaced
  vertical and horizontal circles, it follows that the total signed
  area of any periodic domain is zero. Clearly, a non-zero region with
  this property must have both positive and negative local
  multiplicities.
\end{proof}

For the case of links, we number our dots $\{w_{i,j}\}_{(i,j)\in S}$
and $\{z_{i,j}\}_{(i,j)\in S}$ where $S$ is the index set consisting
of $(i,j)$ with $i=1,...,\ell$ and $j=1,...,n_i$, and the dots
$w_{i,j}$ and $z_{i,j}$ lie on the $i^{th}$ component of $\orL$.

In this case, the Alexander grading is an $\ell$-tuple of integers.
It is uniquely characterized by the property that for $\x,\y\in X$,
the $i^{th}$ component of the Alexander grading is the winding number
of $\gamma_{\x,\y}$ about the sum of the black dots in
$\{z_{i,j}\}_{(i,j)\in S}$ minus its winding number around
$\{w_{i,j}\}_{(i,j)\in S}$. Again, this can be more succinctly
recorded by placing (minus one times) a vector of winding numbers at
each vertex, and defining $A(\x)$ as the sum of these local
contributions at each intersection point in $\x$. This can then be
renormalized to be symmetric.

Once again, we have the chain complex $C(\Gamma)$ as defined in the
introduction, which now inherits an $\ell$-tuple of Alexander gradings
and a single Maslov grading.

There is also a refinement, $C^-(\Gamma)$, which is freely generated
by $X$ over $\Field[\{U_{i,j}\}_{(i,j)\in S}]$, endowed with the
differential
$$\partial^- \x = \sum_{\y\in X} \sum_{r\in R_{\x,\y}}
\left\{\begin{array}{ll} 1 & {\text{if $P_{\x}(r)+P_{\y}(r)=1$}} \\ 0 & {\text{otherwise}} 
\end{array}
\right\}
\cm
\left(\prod_{(i,j)\in S} U_{i,j}^{n_{w_{i,j}}(r)}\right).$$

\begin{theorem}
  There are multi-graded identifications \[
  H_*(C(\Gamma),\partial)\cong {\widehat{HFL}}(\orL)\otimes
  \bigotimes_{i=1}^\ell V_i^{\otimes (n_i-1) }, \] where $V_i$ is the
  two-dimensional vector space spanned by two generators, one in zero
  Maslov and Alexander multigradings, and the other in Maslov grading
  negative one and Alexander multi-grading corresponding to minus the
  $i^{th}$ basis vector.  More generally, the multi-filtered chain
  homotopy type of $\CFLm(S^3,\orL)$ is identified with the
  multi-filtered chain homotopy type of $C^-(\Gamma)$.
\end{theorem}

\begin{proof}
        Both follow from the proof of
        Theorem~\ref{thm:KnotFloerHomology}, combined with
        Proposition~\ref{prop:ExtraBasepointsLink}.
\end{proof}

\section {Examples}
\label{sec:Examples}

We give a few elementary illustrations of our results. 

\subsection{Hopf link}
\label{subsec:Hopf}

\begin{figure}
\begin{center}
\mbox{\vbox{\epsfbox{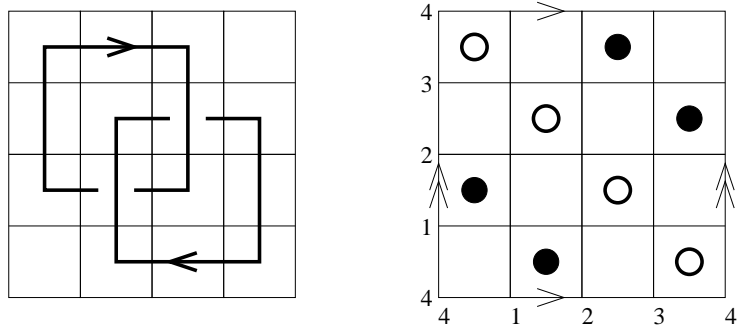}}}
\end{center}
\caption {{\bf Grid presentation of the Hopf link.}}
\label {Figure:hopf}
\end{figure}

Consider the grid
presentation $\Gamma$ for the Hopf link with grid number $n=4,$ shown in
Figure~\ref{Figure:hopf}.

The generators of the chain complex $C(\Gamma)$ are in one-to-one 
correspondence with permutations $\sigma$ of the set $\{1,2,3,4\}.$
For conciseness, we write the generator consisting of the intersections of the 
$i$th horizontal circle with the $\sigma(i)$th vertical circle as 
$\bigl(\sigma(1)\sigma(2)\sigma(3)\sigma(4) \bigr).$

There are eight empty squares in the grid. Each of them produces differentials 
between generators that differ by a transposition, according to the recipe:

$$\begin {array}{ccccccc}
(1\ 2 **) &  & (*\ 2\ 3\ *) &  & (** 3\ 4) &   & (1**\ 4) \\
\downarrow &   & \downarrow &   &
\downarrow  &   & \downarrow \\
(2\ 1 **) &   & (*\ 3\ 2\ *) &  & (** 4\ 3) &   & (4**\ 1) 
\end {array}$$
\medskip

$$\begin {array}{ccccccc}
(3\ 4 **) &   & (*\ 4\ 1\ *) & & (** 1\ 2) & & (3**\ 2) \\
\downarrow &   & \downarrow &   &
\downarrow  &   & \downarrow \\
(4\ 3 **) &  & (*\ 1\ 4\ *) &  & (** 2\ 1) &  & (2**\ 3).
\end {array}$$

\medskip

\begin{figure}
\begin{center}
\mbox{\vbox{\epsfbox{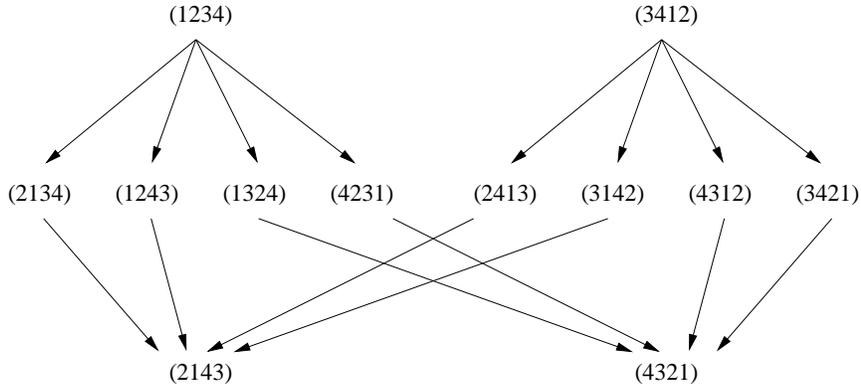}}}
\end{center}
\caption {{\bf Part of the chain complex for the Hopf link.}
This complex appears in Alexander bigrading $(-\frac{1}{2},-\frac{1}{2})$. 
Its homology has rank four.}
\label{Figure:complex}
\end{figure}

The result is that there are sixteen differentials in $(C(\Gamma),
\partial).$ They connect twelve of the $24$ generators, as shown
schematically in Figure~\ref{Figure:complex}. Each of the other twelve
generators is not connected by differentials to any other generators.
Therefore, the homology of our complex has total rank $4+12=16.$

The Alexander bigrading of the generators is computed using the 
adaptation for links of Equation (\ref{eq:Alexanders}). For example,
$$ A(3214) =  \Bigl(\frac{1}{2}, \frac{1}{2} \Bigr), \ \ A(2143) = 
\Bigl( -\frac{1}{2},-\frac{1}{2} \Bigr).$$ 

To compute the Maslov gradings, we start with the canonical generator 
$(2143),$ which has $M=-3.$ Then, we relate each of the other generators 
to the canonical one by a sequence of transpositions. Whenever two 
generators $\mathbf{x}$ and $\mathbf{y}$ differ by a transposition, if a 
two-chain $D$ has boundary $\gamma_{\mathbf{x}, \mathbf{y}},$ then $D$ 
consists of two points and a rectangle, and it is straightforward to apply 
Equation~\eqref{eq:RelativeMaslovGrading}. For example, 
$$M(2134)=M(2143)+1=-2, \ \ M(2314)=M(2134)+1=-1, \ \text{etc.}$$

The result is that
$$ H_*(C(\Gamma), \partial) = \Bigl(V_1^{\otimes 2}\otimes V_2^{\otimes 
2}\Bigr) \Bigl[\frac{1}{2},\frac{1}{2}\Bigr].$$

Here, the notation $[i,j]$ denotes an upward shift in Alexander 
bigrading, i.e. if $V$ is a bigraded vector space, then $(V[i,j])_{x,y} = 
V_{x-i, y-j}.$

The link Floer homology of $\vec H$ is $(V_1\otimes 
V_2)[\frac{1}{2},\frac{1}{2}],$ cf. \cite{Links}. 
This confirms that $$ 
H_*(C(\Gamma), \partial) = \HFLa(\vec H) \otimes V_1\otimes V_2.$$ 

\subsection{The trefoil}

Consider the grid presentation of the trefoil knot shown in
Section~\ref{sec:Introduction}.  There are, of course, $120$
generators of the chain complex.  A quick glance at 
Figure~\ref{Figure:tref}
reveals 15 rectangles containing no black or white dot: fifteen $1\times
1$, five $2\times 1$, and five $1\times 2$. Each rectangle gives rise
to $3!=6$ different differentials. With a little computer assistance
or a great deal of patience, one finds that the homology of this
complex has rank $48$. Indeed, with the conventions used in
Subsection~\ref{subsec:Hopf}, one finds that the generators correspond
to permutations (listed in increasing Alexander grading):
\[
\begin{array}{llllllll}
 (23451) & (13452) & (23415) & (23541) & (24351) & (32451) & (13542) & (14352) \\
 (24153) & (24315) & (25413) & (32415) & (32541) & (35421) & (42351) & (43152) \\
 (43521) & (15342) & (15423) & (25143) & (31542) & (32514) & (35241) & (41352) \\
 (42153) & (42315) & (43125) & (45321) & (52341) & (52413) & (54312) & (15243) \\
 (15324) & (31524) & (41325) & (42135) & (51342) & (51423) & (52143) & (52314) \\
 (54132) & (54213) & (15234) & (41235) & (51243) & (51324) & (52134) & (51234).
\end{array}
\]
These generators have Alexander gradings between $-5$ and $1$; there
are $1, 5, 11, 14, 11, 5, 1$ generators in gradings $-5,-4,-3,-2,-1,0,1,$ 
respectively. The knot Floer homology group for the left-handed trefoil 
$T$ is non-trivial in only three Alexander-Maslov bigradings $(-1,0)$, 
$(0,1)$, and $(1,2)$, and it has rank one in these three bigradings.  
Considering Maslov gradings as well, one immediately verifies that
$$H_*(C(\Gamma))\cong \HFKa(T)\otimes V^{\otimes 4}.$$

\bibliographystyle{plain}
\bibliography{biblio}

\end{document}